\documentclass[a4paper,12pt]{article}
\usepackage[utf8]{inputenc}
\usepackage[]{amsmath}
\usepackage[]{amssymb}
\usepackage[]{amsthm}
\author{Magdalena Zielenkiewicz}
\title{Integrating Schur polynomials using iterated residues at infinity}
\date{}
\def\T{\mathbb{T}}
\def\C{\mathbb{C}}

\def\Q{\mathcal{Q}}
\def\R{\mathcal{R}}
\def\E{\mathbb{E}}
\def\B{\mathbb{B}}
\def\X{X}

\newtheorem{thm}{Theorem}

\newtheorem{formula}{Formula}
\theoremstyle{definition}

\begin{document}
\maketitle

\begin{abstract}
In this paper we show examples of computations achieved using the formulas of our previous paper, which express the push-forwards in equivariant cohomology as iterated residues at infinity. We consider the equivariant cohomology of the complex Lagrangian Grassmannian $LG(n)$ and the orthogonal Grassmannian with the action of the maximal torus. In particular, we show how to obtain some well-known results due to P. Pragacz and J. Ratajski on integrals of Schur polynomials over the Lagrangian Grassmannian $LG(n)$ and the orthogonal Grassmannian $OG(n)$.
\end{abstract}

\section{Introduction}

In our previous paper we have presented a new approach to push-forwards in equivariant cohomology of homogeneous spaces of classical Lie groups. For the homogeneous space $G/P$ of a classical semisimple Lie group $G$ and its maximal parabolic subgroup $P$ with an action of a maximal torus $T$ in $G$ we expressed the Gysin homomorphism $\pi_*: H_{T}^*(G/P) \to H_T^*(pt)$ associated to a constant map $\pi: G/P \to pt$ in the form of an iterated residue at infinity of a certain complex function. The residue formula depended only on the combinatorial properties of the homogeneous space involved, in particular on the Weyl groups of $G$ and $P$ acting on the set of roots of $G$. For details see our previous paper \cite{mz}. The formulas were obtained using localization techniques. In our recent preprint \cite{mz2} we showed how to obtain them in the context of symplectic reductions, using the Jeffrey-Kirwan nonabelian localization theorem. \newline

Residue type formulas for push-forward in equivariant cohomology provide a useful computational tool, reducing the question of computing the push-forward to computing one coefficient of a Laurent series expansion. In this paper we show an example of such a calculation, reproving a classical result of P. Pragacz and J. Ratajski on push-forwards of Schur polynomials in the equivariant cohomology of the Lagrangian Grassmannian \cite{pragacz}. The original proof involves an acute geometrical investigation, whereas the proof presented here is a purely algebraic straightforward computation. 

\section{Preliminaries}

\subsection{Localization in equivariant cohomology}
Assume $X$ is a compact manifold equipped with a torus action and consider the $\T$-equivariant cohomology $H^*_{\T} X $ of $X$. Let $\pi: X \to pt$ be the projection to a point, and consider the equivariant Gysin homomorphism associated with it
\[\pi_* : H^*_{\T}(X) \to H^*_{\T}(pt).\]  
The Gysin homomorphism for the projection $\pi$ as above is often called the equivariant push-forward to a point or integration along fibers, and is denoted by $a \mapsto \int_X a$. \newline

The localization theorem due to D. Quillen asserts that after localizing with respect to the multiplicative set generated by the nonzero characters of the action the $\T$-equivariant cohomology becomes isomorphic to the cohomology of the fixed point set. 

\begin{thm}[Quillen, \cite{quillen}] Let $\X$ be a compact $\T$-space. The inclusion  $i: \X^\T \hookrightarrow \X$ induces an isomorphism
\[ i^*: H_{\T}^*(\X)[(\T^{\#} \setminus\{ 0\})^{-1}] \stackrel{\simeq}{\longrightarrow} H_{\T}^*(\X^\T)[(\T^{\#} \setminus\{ 0\})^{-1}]\]
after localizing with respect to the multiplicative system consisting of finite products of elements $c_1^{\T}(L_{\chi})$, for $\chi \in \T^{\#}\setminus \{0\}$. The bundle $L_{\chi}: \E\T \times^{\T} \C_{\chi} \to \B\T$ is a line bundle determined by the one dimensional representation $\C_{\chi}$ of $\T$ associated to the character $\chi$.
\end{thm}

If the space $X$ is a compact manifold and moreover the set of fixed points is finite, then one of the consequences of the above isomorphism is the following integration formula, due to M. Atiyah and R. Bott and indepentently to N. Berline and M. Vergne.

\begin{thm}[Atiyah, Bott \cite{ab}, Berline, Vergne, \cite{bv}] Suppose $\X$ is a compact manifold with an action of a torus $\T$ such that $\# \X^\T < \infty$. For $x \in H_{\T}^*(\X)$ we have
\[ \int_{\X} x  = \sum_{p \in \X^\T} \frac{i_p^*x}{e_p},\]
where $e_p$ is the equivariant Euler class of the tangent bundle at the fixed point $p$, and $i_p^*: p \hookrightarrow \X$ is induced by the inclusion of the fixed point $p$.
\end{thm}

Using this formula one can express the push-forward for certain homogeneous spaces of Lie groups in the form of an iterated residue at infinity. Examples include classical, Lagrangian and orthogonal Grassmnnians. Within this paper we will mostly be interested in the Lagrangian Grassmannian $LG(n)$. The case of orthogonal Grassmannians is analogous. An example is given in chapter \ref{schur}.

\subsection{Push-forward formulas}

Let $LG(n)$ be the Lagrangian Grassmannian, parametrizing the maximal isotropic subspaces of $\C^{2n}$ equipped with the standard symplectic form, and consider the action of its maximal torus. The fixed points of the torus action can be parametrized using the subsets $I \subseteq \{ 1,\dots,n\}$:
\[p_I = Span\{ q_i, p_j: i \in I, j \notin I \},\]
where $q_1, \dots, q_n, p_n, \dots, p_1$ are the coordinates on $\C^{2n}$. \newline

The weights of the torus action on the tangent space are equal to
\[\{ \pm t_i \pm t_j: 1\leq i<j \leq n\}\cup\{\pm 2t_i: i=1,...,n\}, \] 
where the + sign appears whenever $i,j \in I$.
In this case, the Atiyah-Bott-Berline-Vergne formula \cite{ab} gives:
\[ \int_{LG(n)} \phi(\mathcal{R}) = \sum_{I}\frac{V(t_i, -t_j: i \in I, j \notin I)}{\prod_{i,j=1}^n (\pm2t_i)( \pm t_i \pm t_j)}, \]
where $\phi(\mathcal{R})$ is a characteristic class of the tautological bundle, which at the fixed points of the action is given by a $W$-symmetric polynomial $V$, where $W$ denotes the Weyl group of $LG(n)$. The right-hand side can be expressed as a residue at infinity as follows:

\begin{formula}
\[\int_{LG(n)} \phi(\mathcal{R}) =  Res_{\mathbf{z} = \infty} \frac{V(z_1,...,z_n)\prod_{i<j}(z_j - z_i)}{\prod_{i=1}^n(t_i - z_i)(t_i + z_i)\prod_{i<j}(t_i + t_j)(t_j - t_i)}.\] \label{formula-lagrangian}
\end{formula}
The residue at $\mathbf{z} = \infty$ denotes the iterated residue at infinity in $\mathbf{z} = (z_1,\dots,z_n)$, that is $Res_{\mathbf{z} = \infty}f(\mathbf{z}) = Res_{z_1 = \infty}\dots Res_{z_n = \infty} f(z_1,\dots, z_n)$. \newline

The proof can be found in \cite{mz}. One can similarly describe the push-forward in the equivariant cohomology of the orthogonal Grassmannians $OG(n,2n)$ and $OG(n, 2n+1)$, for which the residue-type push-forward formulas are the following

\[\int_{OG(n,2n)} \phi(\mathcal{R}) =  Res_{\mathbf{z} = \infty} \frac{V(z_1,...,z_n)\prod_{i<j}(z_j - z_i)2^{n-1} z_1...z_n}{\prod_{i=1}^n(t_i - z_i)(t_i + z_i)\prod_{i<j}(t_i + t_j)(t_j - t_i)}\]
\[\int_{OG(n,2n+1)} \phi(\mathcal{R}) =  Res_{\mathbf{z} = \infty} \frac{V(z_1,...,z_n)\prod_{i<j}(z_j - z_i)2^n}{\prod_{i=1}^n(t_i - z_i)(t_i + z_i)\prod_{i<j}(t_i + t_j)(t_j - t_i)}.\]

These formulas are very similar to the one for the Lagrangian Grassmannian and it is very easy to rephrase the results obtained for $LG(n)$ to obtain the analogous formulas for the orthogonal Grassmannians. Hence,  all the computations will be performed in the Lagrangian case from which the orthogonal cases follow directly.

\subsection{Schur polynomials}

A non-increasing sequence $\lambda = (\lambda_1 \geq \lambda_2 \geq \dots)$ of non-negative integers with only finitely many non-zero terms is called a partition of length $n = \# \{i: \lambda_i > 0\}$ and weight $|\lambda| = \sum_i \lambda_i$. To every partition $\lambda$ we can assign a certain polynomial, called the Schur polynomial, defined as follows.

\[ s_{(\lambda_1,\dots ,\lambda_n)}(z_1,\dots,z_n) =  \frac{\det \left[\begin{array}{c c c c}
z_1^{\lambda_1 +n-1} & z_2^{\lambda_1 +n-1} & \dots &  z_n^{\lambda_1 +n-1}\\ 
z_1^{\lambda_2 +n-2} & z_2^{\lambda_2 +n-2} & \dots &  z_n^{\lambda_2 +n-2}\\ 
\vdots & \vdots & \ddots & \vdots \\ 
z_1^{\lambda_n} & z_2^{\lambda_n} & \dots & z_n^{\lambda_n} \end{array}\right]} {\prod_{i<j}(z_j - z_i)} \]

Schur polynomials were first introduced by Jacobi in $\cite{jacobi}$. These polynomials form a basis of the space of symmetric polynomials with integer coefficients which is particularly useful for applications in representation theory and in the study of the geometry of complex Grassmannians. In particular, Schur polynomials are the characters of the finite-dimensional irreducible representations of $GL_n$. More importantly from our point of view, Schur polynomials represent the fundamental cohomology classes of the Schubert subvarieties $Y_{\lambda}$ of complex Grassmannians. Let $Gr_d(V)$ denote the Grassmannian of $d$-dimensional subspaces of an $n$-dimensional vector space $V$ and let $V_{\bullet}: V_1 \subset V_2 \subset \dots \subset V_n=V$ be a fixed complete flag in $V$. For any partition $\lambda$ with at most $d$ parts such that $n-d \geq \lambda_1 \geq \dots \geq \lambda_d \geq 0$ we define the Schubert variety

\[Y_{\lambda} = \{X \in Gr_d(V): \dim(X \cap V_{n-d + i - \lambda_i}) \geq i \textrm{ for } i=1,\dots,d \}.\]
The subvariety $Y_{\lambda}$ has codimension $|\lambda|$ and its class in $H^*(Gr_d(V))$ is given by the Schur polynomial
\[ [ Y_{\lambda}] = s_{\lambda}(x_1,\dots,x_d),\]
where $x_1,\dots,x_d$ are the Chern roots of the tautological bundle. The classes $[Y_{\lambda}]$ indexed by partitions $\lambda$ as described above form an additive basis of the cohomology ring of $Gr_d(V)$. \newline

There are several other useful ways of defining the Schur polynomials, equivalent to the definition above. For example, one can define the Schur polynomials using Young tableaux. A Young diagram associated to a partition $\lambda$ is a collection of boxes arranged in left-justified rows, with $i$-th row having length $\lambda_i$. By filling the boxes with natural numbers $1,\dots,d$ in such a way that the entries weakly increase along rows and strictly increase down each column we obtain a so-called semistandard Young tableau. The Schur polynomials can be then defined as 
\[s_{\lambda}(x_1\dots,x_d) = \sum_{T} x^T = \sum_T x_1^{t_1}\dots x_d^{t_d},\]
where the summation is over the semistandard Yound tableaux $T$ of shape $\lambda$ and $t_i$ is the number of times $i$ appears in $T$. This approach can be found in \cite{fulton}. In our computations we will be using the Jacobi definition.

\section{Push-forward of Schur polynomials}
Recall that for the Lagrangian Grassmnnian we have
\[\int_{LG(n)} \phi(\mathcal{R}) =  Res_{\mathbf{z} = \infty} \frac{V(z_1,...,z_n)\prod_{i<j}(z_j - z_i)}{\prod_{i=1}^n(t_i - z_i)(t_i + z_i)\prod_{i<j}(t_i + t_j)(t_j - t_i)}.\]

We will use this formula to prove the following result \cite{pragacz}:
Let $\lambda$ be a partition of length $< n+1$ and let $\rho(n)$ be the partition $(n,n-1,\ldots,1)$. Let $s_{\lambda}$ denote the Schur polynomial corresponding to the partition $\lambda$. Furthermore, let $\omega: LG(n)\to pt$ and let
 \[ \omega_*:H_{\T}^*(LG(n))\to H_{\T}^*(pt)\] be the integration over the fiber of $\omega$. We will show the following result.
\begin{thm}
Let $\Q$ be the tautological quotient rank n bundle on $LG(n)$. Then the Schur polynomial $s_{\lambda}(\Q)$ has a nonzero image under $\omega_*$ only if $\lambda = 2\mu+\rho(n)$ for some partition $\mu$. In terms of Young diagrams, this means that the diagram corresponding to $\lambda$ contains the diagram of the standard partition $\rho(n)$, and has additionally an even number of boxes added in each row. In this case, the image is:
\[\omega_* s_{\lambda}(\Q)=s_{\mu}^{[2]}(\C^{2n}),\]
where $s_{\mu}^{[2]}$ is obtained from $s_{\mu}$  by replacing each $e_i$ in the presentation of $s_{\mu}$ as a polynomial in elementary symmetric functions $e_i$ by $(-1)^i c_{2i}(\C^{2n})$. In other words, $s_{\mu}^{[2]}(t_1, \dots, t_n)= s_{\mu}(t_1^2, \dots, t_n^2)$. \label{pr}
\end{thm}

Note that even though formally we only state the theorem for the push-forward in equivariant cohomology it actually provides us with a tool to study push-forwards in a more general setting. The push-forward in $\T$-equivariant cohomology we describe is the push-forward for the fiber bundle $\E\T \times^{\T} LG(n) \to \B\T$ with fiber $LG(n)$. Since $\B\T$ is the universal base for $\T$-principal bundles we can extend our result to any $\T$-principal bundle $E \to X$ with fiber $LG(n)$, because any such bundle is a pullback of the bundle $\E\T \times^{\T} LG(n) \to \B\T$ (for a precise statment and proof see \cite{tu}). Moreover, the $\T$-principal bundles are precisely the split bundles, thus we can deduce the general result from the splitting principle. We reprove the well-known result of P. Pragacz and J. Ratajski on push-forwards of Schur polynomials to show the computational advantages of our residue-type formulas - the proof we present here is short and purely algebraic, in contrast with the insightful geometric analysis in \cite{pragacz}. 

\subsection{The case $n=2$}
Before we prove it, let us look at an example for $n=2$. In this case formula \ref{formula-lagrangian} has the form:
\[\int_{LG(2)} s_{(\lambda_1,\lambda_2)}(\R) = \]
\[= Res_{z_1=z_2 = \infty} \frac{s_{(\lambda_1, \lambda_2)}(z_1,z_2)(z_2 - z_1)}{(t_1 - z_1)(t_1 + z_1)(t_2 - z_2)(t_2 + z_2)(t_1 + t_2)(t_2 - t_1)}= \]
\[=\frac{1}{t_2^2 - t_1^2} Res_{z_1=z_2 = \infty} \frac{\det \left[\begin{array}{c c}
z_1^{\lambda_1 +1} & z_2^{\lambda_1 +1} \\ 
z_1^{\lambda_2} & z_2^{\lambda_2} \end{array}\right] \frac{(z_2-z_1)}{(z_2 - z_1)} }{(t_1 - z_1)(t_1 + z_1)(t_2 - z_2)(t_2 + z_2)} = \frac{1}{t_2^2 - t_1^2} \cdot  \bigstar \]

To compute the residue we use the fact that $Res_{z=\infty} f(z) = -\frac{1}{z^2} Res_{z=0} f(\frac{1}{z})$, so

\[\bigstar = \frac{1}{z_1^2 z_2^2} Res_{z_1=z_2 = 0} \frac{\det \left[\begin{array}{c c}
z_1^{-(\lambda_1 +1)} & z_2^{-(\lambda_1 +1)} \\ 
z_1^{-\lambda_2} & z_2^{-\lambda_2} \end{array}\right]}{(t_1 - z_1^{-1})(t_1 + z_1^{-1})(t_2 - z_2^{-1})(t_2 + z_2^{-1})}  = \]
\[=Res_{z_1=z_2 = 0} \frac{\det \left[\begin{array}{c c}
z_1^{-(\lambda_1 +1)} & z_2^{-(\lambda_1 +1)} \\ 
z_1^{-\lambda_2} & z_2^{-\lambda_2} \end{array}\right]}{(t_1 z_1- 1)(t_1 z_1 + 1)(t_2 z_2 - 1)(t_2 z_2 + 1)} \]

Note that the residue at zero of a function is equal to the coefficient corresponding to $z^{-1}$ in the Laurent series expansion. Thus, in order to compute $\bigstar$ we need to expand as a power series the function 
\[\frac{1}{(t_1 z_1- 1)(t_1 z_1 + 1)(t_2 z_2 - 1)(t_2 z_2 + 1)} = \frac{1}{((t_1 z_1)^2- 1)((t_2 z_2)^2 - 1)},\]
 multiply by the determinant coming from the Schur polynomial, and take the coefficient corresponding to $z_1^{-1} z_2^{-1}$. \newline

Proceeding in the steps described above, we have:
\[\frac{1}{((t_1 z_1)^2- 1)((t_2 z_2)^2 - 1)} = \sum_{i=0}^{\infty} \sum_{j=0}^{\infty} t_1^{2 i} t_2^{ 2 j} z_1^{2 i} z_2^{2 j} = \bigstar \bigstar.\]

Now we expand the determinant: 
\[\det \left[\begin{array}{c c}
z_1^{-(\lambda_1 +1)} & z_2^{-(\lambda_1 +1)} \\ 
z_1^{-\lambda_2} & z_2^{-\lambda_2} \end{array}\right] = z_1^{-(\lambda_1 + 1)} z_2^{-\lambda_2} - z_2^{-(\lambda_1 + 1)} z_1^{-\lambda_2} = A + B\]

From this we can easily deduce the first claim of the Theorem \ref{pr}: \newline

The coefficients in the series $\bigstar \bigstar$ are even, so if $\lambda_1 + 1$ is odd or $\lambda_2$ is even, then the coefficient at $z_1^{-1} z_2^{-1}$ in the result is zero, because it is zero both in $A \cdot \bigstar \bigstar$ and in $B \cdot \bigstar \bigstar$.  Also, if $\lambda_1 < 2$ or $\lambda_2 < 1$ then the result is zero, since $(A+B)\cdot \bigstar \bigstar$ is nonsingular at $0$. This shows that the expression $\bigstar$ can only be nonzero if $\lambda_1$ is even and greater than $2$, $\lambda_2$ is odd and greater than $1$, so $\lambda = (2,1)+2 \mu $. \newline
In this case, $\bigstar$ splits into a sum of two contributions coming from $A$ and $B$:

Contribution from $A$:
\[X_A=A \cdot \bigstar \bigstar = z_1^{-(\lambda_1 + 1)} z_2^{-\lambda_2} \sum_{i=0}^{\infty} \sum_{j=0}^{\infty} t_1^{2 i} t_2^{ 2 j} z_1^{2 i} z_2^{2 j} = \sum_{i=0}^{\infty} \sum_{j=0}^{\infty} t_1^{2 i} t_2^{ 2 j} z_1^{2 i - (\lambda_1 + 1)} z_2^{2 j - \lambda_2}\]
The coefficient at $z_1^{-1 z_2^{-1}}$ is equal to $t_1^{\lambda_1} t_2^{\lambda_2 -1}$. \newline

Similarly, the contribution from $B$ is:
\[X_B=B \cdot \bigstar \bigstar = z_2^{-(\lambda_1 + 1)} z_1^{-\lambda_2} \sum_{i=0}^{\infty} \sum_{j=0}^{\infty} t_1^{2 i} t_2^{ 2 j} z_1^{2 i} z_2^{2 j} = \sum_{i=0}^{\infty} \sum_{j=0}^{\infty} t_1^{2 i} t_2^{ 2 j} z_1^{2 i - \lambda_2 } z_2^{2 j - (\lambda_1+1)}\]
The coefficient at $z_1^{-1} z_2^{-1}$ is equal to $t_1^{\lambda_2 - 1} t_2^{\lambda_1}$. \newline

Finally, \[\bigstar = X_A - X_B = t_1^{\lambda_1} t_2^{\lambda_2 -1} - t_1^{\lambda_2 - 1} t_2^{\lambda_1} = \det \left[\begin{array}{c c}
t_1^{\lambda_1} & t_2^{\lambda_1} \\ 
t_1^{\lambda_2 - 1} & t_2^{\lambda_2 - 1} \end{array}\right], \]
and so the push-forward is
\[\omega_* s_{\lambda}(\R) = \frac{1}{t_2^2 - t_1^2} \cdot  \bigstar = \frac{\det \left[\begin{array}{c c}
t_1^{\lambda_1} & t_2^{\lambda_1} \\ 
t_1^{\lambda_2 - 1} & t_2^{\lambda_2 - 1} \end{array}\right]}{t_2^2 - t_1^2} =\]
\[= \frac{\det \left[\begin{array}{c c}
(t_1^2)^{\lambda_1 / 2 } & (t_2^2)^{\lambda_1 / 2} \\ 
(t_1^2)^{(\lambda_2 - 1)/2} & (t_2^2)^{(\lambda_2 - 1)/2} \end{array}\right]}{t_2^2 - t_1^2} 
= \frac{\det \left[\begin{array}{c c}
(t_1^2)^{\mu_1 + 1 } & (t_2^2)^{\mu_1 + 1} \\ 
(t_1^2)^{\mu_2} & (t_2^2)^{\mu_2} \end{array}\right]}{t_2^2 - t_1^2} = s_{\mu}(t_1^2, t_2^2),\]
which is exactly what we wanted.

\subsection{Arbitrary $n$}\label{schur}
The case of $n=2$ can be easily generalized to work for any $n$. We need to proceed in exactly the same steps:
\begin{itemize}
	\item Use the iterated residue formula \ref{formula-lagrangian}. 
	\item Change variables to compute the residue at $0$ and simplify the expression.
	\item Expand the denominator into a power series.
	\item Expand the determinant coming from the Schur polynomial, compute the contributions from the summands, and add them up.
\end{itemize}

Let us follow this procedure. 

\[\int_{LG(n)} s_{(\lambda_1,\dots ,\lambda_n)}(\R) = \]
\[= Res_{z_1=\dots = z_n = \infty} \frac{s_\lambda (z_1,\dots, z_n)\prod_{i<j}(z_j - z_i)}{\prod_{i=1}^n(t_i - z_i)(t_i + z_i)\prod_{i<j}(t_i + t_j)(t_j - t_i)}= \]
\[=\frac{1}{\prod_{i<j}(t_j^2 - t_i^2)} Res_{z_1=\dots=z_n = \infty} \frac{\det \left[\begin{array}{c c c c}
z_1^{\lambda_1 +n-1} & z_2^{\lambda_1 +n-1} & \dots &  z_n^{\lambda_1 +n-1}\\ 
z_1^{\lambda_2 +n-2} & z_2^{\lambda_2 +n-2} & \dots &  z_n^{\lambda_2 +n-2}\\ 
\vdots & \vdots & \vdots & \vdots \\ 
z_1^{\lambda_n} & z_2^{\lambda_n} & \dots & z_n^{\lambda_n} \end{array}\right] }{\prod_{i=1}^n(t_i - z_i)(t_i + z_i)} =\]
\[= \frac{1}{\prod_{i<j}(t_j^2 - t_i^2)} \cdot  \bigstar \]

Now we change variables to compute the residue at $0$:
\[ \bigstar =  \frac{(-1)^n}{z_1^2 \dots z_n^2} Res_{z_1 = \dots = z_n = 0}\frac{\det \left[\begin{array}{c c c c}
z_1^{-(\lambda_1 +n-1)} & z_2^{-(\lambda_1 +n-1)} & \dots &  z_n^{-(\lambda_1 +n-1)}\\ 
z_1^{-(\lambda_2 +n-2)} & z_2^{-(\lambda_2 +n-2)} & \dots &  z_n^{-(\lambda_2 +n-2)}\\ 
\vdots & \vdots & \vdots & \vdots \\ 
z_1^{-\lambda_n} & z_2^{-\lambda_n} & \dots & z_n^{-\lambda_n} \end{array}\right] }{\frac{1}{z_1^2 \dots z_n^2} \prod_{i=1}^n(t_i z_i - 1)(t_i z_i + 1)} \]

The series expansion of the denominator is:
\[\frac{1}{((t_1 z_1)^2- 1) \dots ((t_n z_n)^2 - 1)} = \sum_{i_1=0}^{\infty} \dots \sum_{i_n=0}^{\infty} t_1^{2 i_1} \dots t_n^{ 2 i_n} z_1^{2 i_1} \dots z_n^{2 i_n} = \bigstar \bigstar.\]

Now we expand the determinant as the sum over all permutations:

\[\det \left[\begin{array}{c c c c}
z_1^{-(\lambda_1 +n-1)} & z_2^{-(\lambda_1 +n-1)} & \dots &  z_n^{-(\lambda_1 +n-1)}\\ 
z_1^{-(\lambda_2 +n-2)} & z_2^{-(\lambda_2 +n-2)} & \dots &  z_n^{-(\lambda_2 +n-2)}\\ 
\vdots & \vdots & \vdots & \vdots \\ 
z_1^{-\lambda_n} & z_2^{-\lambda_n} & \dots & z_n^{-\lambda_n} \end{array}\right] = \]
\[=  \sum_{\sigma \in \Sigma_n} (-1)^{ sgn(\sigma)} z_{\sigma(1)}^{-(\lambda_1 + n -1)} z_{\sigma(2)}^{-(\lambda_2 + n - 2)}\dots z_{\sigma(n)}^{-\lambda_n}, \]
and compute the contribution coming from one summand:
\[ (-1)^{ sgn(\sigma)} z_{\sigma(1)}^{-(\lambda_1 + n -1)} z_{\sigma(2)}^{-(\lambda_2 + n - 2)}\dots z_{\sigma(n)}^{-\lambda_n} \cdot \bigstar \bigstar =\]
\[= (-1)^{ sgn(\sigma)} z_{\sigma(1)}^{-(\lambda_1 + n -1)} z_{\sigma(2)}^{-(\lambda_2 + n - 2)}\dots z_{\sigma(n)}^{-\lambda_n} \cdot \sum_{i_1, \dots, i_n=0}^{\infty}  t_1^{2 i_1} \dots t_n^{ 2 i_n} z_1^{2 i_1} \dots z_n^{2 i_n}= \]
\[= (-1)^{ sgn(\sigma)} z_{\sigma(1)}^{-(\lambda_1 + n -1)} z_{\sigma(2)}^{-(\lambda_2 + n - 2)}\dots z_{\sigma(n)}^{-\lambda_n} \cdot \!\!\!\!\! \sum_{i_{\sigma(1)}, \dots, i_{\sigma(n)} =0}^{\infty} \!\!\!\!\! t_{\sigma(1)}^{2 i_{\sigma(1)}} \dots t_{\sigma(n)}^{ 2 i_{\sigma(n)}} z_{\sigma(1)}^{2 i_{\sigma(1)}} \dots z_{\sigma(n)}^{ 2 i_{\sigma(n)}}\]

The coefficient at $z_1^{-1} \dots z_n^{-1}$ is $t_{\sigma(1)}^{2 i_{\sigma(1)}} \dots t_{\sigma(n)}^{ 2 i_{\sigma(n)}}$, where $i_{\sigma(j)}$ must satisfy:

\[ 2 i_{\sigma(j)} - (\lambda_j + n - j) = -1 \textrm{ for } j=1,\dots,n.\]

The solutions are $2 i_{\sigma(j)} =\lambda_j + n - j -1 $. In particular this shows that the coefficient we are looking for is always zero unless $\lambda_j + n - j -1$ is even, i.e.
\[\lambda_j = 2 k_j - n + j + 1 = 2 k_j - n + j + 1 + (n+j-1) - (n+j-1) = 2 k_j - 2 n + 2 +\rho(n)_j,\]
so $\lambda$ is of the form $\rho(n) + 2\mu$ for some partition $\mu$. \newline

Finally, if $\lambda = \rho(n) + 2\mu$, the sum of all contributions is
\[\bigstar =  \sum_{\sigma \in \Sigma_n} (-1)^{ sgn(\sigma)} t_{\sigma(1)}^{2 i_{\sigma(1)}} \dots t_{\sigma(n)}^{ 2 i_{\sigma(n)}} = \sum_{\sigma \in \Sigma_n} (-1)^{ sgn(\sigma)} t_{\sigma(1)}^{\lambda_1 + n - 2}t_{\sigma(2)}^{\lambda_2 + n - 3} \dots t_{\sigma(n)}^{ \lambda_n -1}= \]
\[= \det \left[\begin{array}{c c c c}
t_1^{\lambda_1 +n-2} & t_2^{\lambda_1 +n-2} & \dots &  t_n^{\lambda_1 +n-2}\\ 
t_1^{\lambda_2 +n-3} & t_2^{\lambda_2 +n-3} & \dots & t_n^{\lambda_2 +n-3}\\ 
\vdots & \vdots & \vdots & \vdots \\ 
t_1^{\lambda_n-1} & t_2^{\lambda_n-1} & \dots & t_n^{\lambda_n-1} \end{array}\right] \]

The resulting expression for the push-forward is then:

\[ \frac{1}{\prod_{i<j}(t_j^2 - t_i^2)} \cdot \det \left[\begin{array}{c c c c}
(t_1^2)^{(\lambda_1 +n-2)/2} & (t_2^2)^{(\lambda_1 +n-2)/2} & \dots &  (t_n^2)^{(\lambda_1 +n-2)/2}\\ 
(t_1^2)^{(\lambda_2 +n-3)/2} & (t_2^2)^{(\lambda_2 +n-3)/2} & \dots & (t_n^2)^{(\lambda_2 +n-3)/2}\\ 
\vdots & \vdots & \vdots & \vdots \\ 
(t_1^2)^{(\lambda_n-1)/2} & (t_2^2)^{(\lambda_n-1)/2} & \dots & (t_n^2)^{(\lambda_n-1)/2} \end{array}\right]=\]
\[= \frac{\det \left[\begin{array}{c c c c}
(t_1^2)^{\mu_1 +n-1} & (t_2^2)^{\mu_1 +n-1)} & \dots &  (t_n^2)^{\mu_1 +n-1}\\ 
(t_1^2)^{\mu_2 +n-2} & (t_2^2)^{\mu_2 +n-2)} & \dots &  (t_n^2)^{\mu_2 +n-2)}\\ 
\vdots & \vdots & \vdots & \vdots \\ 
(t_1^2)^{\mu_n} & (t_2^2)^{\mu_n} & \dots & (t_n^2)^{\mu_n} \end{array}\right]}{\prod_{i<j}(t_j^2 - t_i^2)} = s_{\mu}(t_1^2, \dots, t_n^2),\]
which is the expression we were supposed to prove. \newline

Analogous results can be obtained in the same manner for the orthogonal Grassmannians. For example, for $OG(n, 2n+1)$ one can prove the following formula \cite{pragacz} (using the notation as for the Lagrangian case):
\[\omega_* s_{\lambda}(\Q)=2^n s_{\mu}^{[2]}(\C^{2n}),\]
for a partition $\lambda$ of the form $\lambda = 2\mu + \rho_n$.
This result is a straightforward adaptation of the proof in the Lagrangian case: the formulas differ by the $2^n$ factor only. A slightly different result holds in the even orthogonal case: the push-forward over $OG(n,2n)$ of a Schur polynomial given by the partition $\lambda = 2 \mu + \rho_{n-1}$ equals
\[\omega_* s_{\lambda}(\Q)=2^n s_{\mu}^{[2]}(\C^{2n}).\]
Again, the proof is just an adaptation of the proof in the Lagrangian case: here the formulas differ by the factor $2^{n-1} z_1 \dots z_n$, the product of the variables $z_i$ shifts the degree of each term in the power series expansion by $(1,1,\dots,1)$, hence instead of the partition $\rho_n$, the partition $\rho_{n-1}$ appears.

\section{Push-forward of some other polynomials}

Let us come back to the Lagrangian case. Following the computations in Section \ref{schur} one can easily prove a slightly more general result, extracting the necessary conditions a polynomial defining the characteristic class we wish to push-forward must satisfy in order to make an analogous computation successful. \newline

Consider the characteristic class of the tautological bundle over $LG(n)$ which at the fixed point of the action is given by a $W$-symmetric polynomial $V(z_1, \dots, z_n)$ and assume $V$ can be written as
\[V(z_1,\dots,z_n) = \frac{W(z_1,\dots,z_n)}{\prod_{i < j}(z_j - z_i)}.\] 
Then

\[\int_{LG} V = Res_{\mathbf{z} = \infty} \frac{V(z_1,...,z_n)\prod_{i<j}(z_j - z_i)}{\prod_{i=1}^n(t_i - z_i)(t_i + z_i)\prod_{i<j}(t_i + t_j)(t_j - t_i)} = \]
\[=\frac{1}{\prod_{i<j}(t_j^2 - t_i^2)}\cdot Res_{\mathbf{z} = \infty} \frac{W(z_1,...,z_n)}{\prod_{i=1}^n(t_i - z_i)(t_i + z_i)} = \frac{1}{\prod_{i<j}(t_j^2 - t_i^2)} \cdot \bigstar.\]

Changing the variables $z_i \mapsto \frac{1}{z_i}$ we obtain

\[\bigstar = \frac{(-1)^n}{z_1^2 \dots z_n^2}Res_{\mathbf{z} = \infty} \frac{W(z_1,...,z_n)}{\prod_{i=1}^n(t_i - z_i)(t_i + z_i)} = \]
\[= \frac{(-1)^n}{z_1^2 \dots z_n^2} Res_{\mathbf{z} = 0} \frac{W(z_1^{-1}, \dots, z_n^{-1})}{\prod_{i=1}^n(t_i z_i - 1)(t_i z_i + 1)} = \]
\[= \frac{(-1)^n}{z_1^2 \dots z_n^2} Res_{\mathbf{z} = 0} \frac{W(z_1^{-1}, \dots, z_n^{-1})}{\prod_{i=1}^n(t_i^2 z_i^2 - 1)} = \frac{(-1)^n}{z_1^2 \dots z_n^2} \cdot \heartsuit\]
Expanding the geometric series we get
\[\heartsuit = Res_{\mathbf{z} = 0} (-1)^n W(z_1^{-1}, \dots, z_n^{-1}) \sum_{i_1 = 0}^{\infty} \dots \sum_{i_n = 0}^{\infty} t_1^{2 i_1} \dots t_n^{2 i_n} z_1^{2 i_1} \dots z_n^{2 i_n},\]
hence
\[\bigstar = \frac{(-1)^n}{z_1^2 \dots z_n^2} \cdot \heartsuit = \]
\[= Res_{\mathbf{z} = 0} W(z_1^{-1}, \dots, z_n^{-1}) \sum_{i_1 = 0}^{\infty} \dots \sum_{i_n = 0}^{\infty} t_1^{2 i_1} \dots t_n^{2 i_n} z_1^{2 i_1 - 2} \dots z_n^{2 i_n -2}.\]

Assume the polynomial $W$ is given by  
\[W(z_1, \dots, z_n) = \sum_{k_1, \dots, k_n} a_{k_1,\dots,k_n} z_1^{k_1}\dots z_n^{k_n}.\]
Then by the above computation we have

\[\bigstar = Res_{\mathbf{z} = 0} \sum_{k_1, \dots, k_n} a_{k_1,\dots,k_n} z_1^{-k_1}\dots z_n^{-k_n} \sum_{i_1 = 0}^{\infty} \dots \sum_{i_n = 0}^{\infty} t_1^{2 i_1} \dots t_n^{2 i_n} z_1^{2 i_1 - 2} \dots z_n^{2 i_n -2}= \]
\[ =  Res_{\mathbf{z} = 0} \sum_{k_1, \dots, k_n} a_{k_1,\dots,k_n}  \sum_{i_1 = 0}^{\infty} \dots \sum_{i_n = 0}^{\infty} t_1^{2 i_1} \dots t_n^{2 i_n} z_1^{2 i_1 - 2 - k_1} \dots z_n^{2 i_n -2 - k_n}.\]
The residue at $0$ equals the coefficient at $z_1^{-1}\dots z_n^{-1}$, which is realized by the terms in the above series which are given by conditions:

\[2 i_j -2 - k_j = -1 \textrm{ for } j = 1, \dots, n.\]
Therefore we have

\[  Res_{\mathbf{z} = 0} \sum_{k_1, \dots, k_n} a_{k_1,\dots,k_n}  \sum_{i_1 = 0}^{\infty} \dots \sum_{i_n = 0}^{\infty} t_1^{2 i_1} \dots t_n^{2 i_n} z_1^{2 i_1 - 2 - k_1} \dots z_n^{2 i_n -2 - k_n} = \]
\[= \sum_{k_1, \dots, k_n} a_{k_1,\dots,k_n} t_1^{k_1 + 1} \dots t_n^{k_n + 1}.\]

Finally, the push-forward of $V$ is given by

\[\int_{LG(n)} V  = \frac{1}{\prod_{i<j}(t_j^2 - t_i^2)}\sum_{k_1, \dots, k_n} a_{k_1,\dots,k_n} t_1^{k_1 + 1} \dots t_n^{k_n + 1}.\]

\subsection{Special case}

Assume that $W(z_1,\dots,z_n)$ is a polynomial with all terms having an even degree in each variable (resp. all terms having an odd degree),
\[W(z_1,\dots,z_n) = \sum_{m_1, \dots, m_n} a_{2 m_1,\dots,2 m_n} z_1^{2 m_1} \dots z_n^{2 m_n}\]
(resp. $W(z_1,\dots,z_n) = \sum_{m_1, \dots, m_n} a_{2 m_1 +1,\dots,2 m_n +1} z_1^{2 m_1 +1} \dots z_n^{2 m_n+1}$). Then the push-forward formula reduces to

\[\int_{LG(n)} V  = \frac{1}{\prod_{i<j}(t_j^2 - t_i^2)}\sum_{ m_1, \dots,  m_n} a_{2 m_1,\dots,2 m_n} t_1^{2m_1 + 1} \dots t_n^{2m_n + 1} = \]
\[ = \frac{1}{\prod_{i<j}(t_j^2 - t_i^2)} t_1 \dots t_n \sum_{ m_1, \dots,  m_n} a_{2 m_1,\dots,2 m_n} t_1^{2m_1} \dots t_n^{2m_n } = \]
\[= t_1 \dots t_n \frac{W(t_1^2,\dots, t_n^2)}{\prod_{i<j}(t_j^2 - t_i^2)} = \frac{ t_1 \dots t_n}{{\prod_{i<j}(t_j^2 - t_i^2)}} \tilde{W}(t_1^2, \dots, t_n^2),\]
where the polynomial $\tilde{W}$ is obtained from $W$ by shifting the coefficients: that is, if $W(z_1\dots,z_n) = \sum_{k_1, \dots, k_n} a_{k_1,\dots,k_n} z_1^{k_1}\dots z_n^{k_n}$, then
\[\tilde{W}(z_1\dots,z_n) = \sum_{k_1, \dots, k_n} a_{k_1,\dots,k_n} z_1^{2 k_1}\dots z_n^{2 k_n}.\]

Analogously, for a polynomial $W$ with only odd degree terms in each variable we have
\[\int_{LG(n)} V =  \frac{ t_1^2 \dots t_n^2}{{\prod_{i<j}(t_j^2 - t_i^2)}} \tilde{W}(t_1^2, \dots, t_n^2).\]

\end{document}